\documentclass{article}

\usepackage{amsmath}
\usepackage[francais]{babel}
\usepackage{amssymb}
\usepackage[applemac]{inputenc}
\usepackage{mathrsfs}

\def\noi {\noindent}
\def\Q{\mathop{\mathbb Q}}
\def\F{\mathop{\mathbb F}}

\def\R{\mathop{\mathbb R}}
\def\deg{\mathop{\rm deg}}
\def\dim{\mathop{\rm dim}}
\def\Sym{\mathop{\rm Sym}}
\def\Alt{\mathop{\rm Alt}}

\def\sgn{\mathop{\rm sgn}}
\def\tr{\mathop{\rm tr}}
\def\rg{\mathop{\rm rg}}

\def\GL{\mathop{\rm GL}}

\def\Aut{\mathop{\rm Aut}}
\def\Out{\mathop{\rm Out}}

\def\cc{\mathscr{C}}

\def\sm{\raisebox{1.5pt}{~\rule{5pt}{1.3pt}~}}

\begin{document}
 
  
 \medskip
 \centerline{\bf   Groupes de Coxeter finis : centralisateurs d'involutions} 
 
 \medskip
 \centerline{Springer LNM {\bf2313} (2023), 279-300}
      
  \bigskip

   \smallskip
   
   \hspace{33mm} {\it À Catriona Byrne, en souvenir d'une vieille amitié}
   
   \vspace{5mm} 
   
\noi {\bf Introduction}.

\smallskip Soit $G$  un groupe de Coxeter fini, soit $u$  une involution de $G$ et soit $G_u$ le centralisateur de $u$ dans $G$. Dans certains cas, par exemple quand $u$ est une réflexion, le groupe $G_u$ est engendré par des réflexions de $G$; en particulier, c'est un groupe de Coxeter.  Il n'en est pas de même en général, mais c'est ``presque'' le cas. Notre but est de préciser cet énoncé (cf. th.1.1  ci-dessous),
 et de décrire explicitement certains invariants de $G_u$, pour chaque type $A_n, B_n, ..., I_2(m)$. 
 
  \bigskip
  
  \noi {\bf 1. Premiers énoncés}.
  
  \smallskip
    Rappelons d'abord quelques notations et quelques définitions.
   
 \smallskip  
 
   Soit $V$ un $\R$-espace vectoriel de dimension finie, et soit $G$ un sous-groupe fini de $\GL(V)$ engendré par des réflexions, i.e. par des éléments d'ordre 2 fixant un hyperplan. Nous dirons, comme dans [Se 22], 1.1, que le couple $(V,G)$ est un {\it couple de Coxeter}. Sauf mention expresse du contraire (cf. §4 ou 5), on supposera que $V$ est {\it réduit}, c'est-à-dire ne contient aucun élément  $\neq 0$ fixé par $G$; cela revient à demander que $\dim V$ est égal au rang $\rg(G)$ de $G$.
  
  \smallskip
  
  Soit $H$ un sous-groupe de $G$. On dit que $H$ est un $\cc$-{\it sous-groupe} de $G$ s'il est engendré par des réflexions, i.e. si $(V,H)$ est un couple de Coxeter. On dit que $H$ est {\it parabolique} s'il existe une partie $X$ de $V$ telle que $H$ soit l'ensemble des  éléments de $G$ qui fixent $X$; on sait que cela entraîne que $H$ est un $\cc$-sous-groupe, cf. [Se 22], 1.5.
  
  \smallskip
    Soit  $u$ une involution de $G$, autrement dit un élément de $G$ de carré $1$. On 
     note $V_u^{^+}$ le sous-espace vectoriel de $V$ fixé par $u$, et $ V_u^{^-}$ celui fixé par $-u$. On a $V = V_u^{^+} \oplus V_u^{^-}$.     Le {\it degré} de $u$
    est défini par $\deg(u) = \dim V_u^{^-}$; c'est la multiplicité de $-1$ comme valeur propre de $u$; on le note souvent $d$.
        
  \medskip     
        
    \noi{\bf Théorème 1.1}. {\it Le centralisateur  $G_u$ de $u$ dans $G$ est engendré par des involutions de degré $1$ et $2$.}
    
    \smallskip
    Dans le cas particulier où $u$ est une involution de degré maximal, c'est le cor.3.18 de [Se 22].
  
  \medskip  
    
      Soit $G^1_u$ le sous-groupe de $G_u$ engendré par les éléments de $G_u$ qui sont des réflexions de $G$. C'est le plus grand $\cc$-sous-groupe de $G_u$. Il est normal dans $G_u$. Notons $\Gamma_u$
      le quotient $G_u/G^1_u$; ce groupe précise dans quelle mesure $G_u$ n'est pas
     engendré par des réflexions. Le th.1.1 équivaut à dire que $\Gamma_u$ est
     engendré par les images des involutions de degré 2 de $G_u$.
      
      \medskip
      
      \noi {\it Réduction au cas irréductible}. 
      
      Supposons que $V$ soit réduit. On sait (cf. [Bo  68], V3.7) que
      $V$ se décompose de façon unique en une somme directe $V = \oplus V_i$
      de représentations irréductibles non triviales de $G$, et que $G= \prod G_i$,
      avec $G_i \subset \GL(V_i)$. Les $G_i$ sont les {\it composantes 
      irréductibles} de $G$. On a $G_u =  \prod G_{u_i}$, où les $u_i$ sont les
      composantes de $u$. Il y a des décompositions analogues pour  $G^1_u, \Gamma_u$, etc. En particulier, il suffit de démontrer le th.1.1 lorsque $G$ est irréductible, donc de l'un des types $A, B, ..., I$; c'est ce que nous ferons dans la suite.
      
      \medskip

        \noi{\bf Théorème 1.2}. (a) {\it Si $G$ est irréductible non de type $D_n \ (n \geqslant 5)$, le groupe  $\Gamma_u$ est isomorphe à un groupe symétrique.}
        
        (b) {\it Si $G$ est irréductible de type $D_n$, $\Gamma_u$ est isomorphe, soit à un groupe symétrique, soit au produit d'un groupe symétrique par un groupe d'ordre $2$.}
      
  \noi       {\small [Par exemple, quand $G$ est de type $D_5$, il existe une involution $u$ de $G$ telle que $\Gamma_u$ soit abélien élémentaire de type $(2,2)$.]}
  
     \medskip   
        \noi {\it Notation}. Dans le cas (a), si $\Gamma_u \neq 1$, nous noterons $\gamma_u$ l'unique entier $r>1$ tel que $\Gamma_u \simeq \Sym_r$. Lorsque
         $\Gamma_u = 1$, nous écrirons tantôt $\gamma_u = 1$ et tantôt $\gamma_u = 0$, suivant le contexte.

        \medskip
     Les théorèmes 1.1 et 1.2 seront démontrés dans les $\S\S$3-12 par une analyse cas par cas,
     qui donnera la structure de $\Gamma_u$, ainsi que celle des groupes de Coxeter $ G_u^{^+},  \widetilde{G}_u^{^+}, G_u^{^-}, \widetilde{G}_u^{^-}$ définis
    ci-dessous. Nous verrons également que
        l'on peut choisir l'isomorphisme $\Gamma_u \to \Sym_{\gamma_u}$ du théorème 1.2 (a) 
        de telle sorte que toute transposition de $\Sym_{\gamma_u}$ soit l'image d'une involution de degré $2$ de $G_u$; il y a un énoncé analogue
        dans le cas du type $D_n$, cf. §6.   
        
         \bigskip
        
         \noi {\bf 2. Les groupes $G_u^{^+}, G_u^{^-}, \widetilde{G}_u^{^+}, \widetilde{G}_u^{^-}$}.
           
         \smallskip
         
         L'action de $G_u$ sur $V$ respecte la décomposition $V = V_u^{^+} \oplus V_u^{^-}$. On a donc
  $G_u \ \subset  \ \GL(V_u^{^+})\times \GL(V_u^{^-})$, ce qui permet de définir les
  quatre groupes suivants :
         
       \smallskip
       
 \quad     $ G_u^{^+} = G_u \cap \GL(V_u^{^+})$,
      
 \quad      $ G_u^{^-} = G_u \cap \GL(V_u^{^-})$,
       
  \quad    $ \widetilde{G}_u^{^+} = $ image de $G_u$ dans $\GL(V_u^{^+})$ par la première projection,
        
  \quad        $ \widetilde{G}_u^{^-} = $ image de $G_u$ dans $\GL(V_u^{^-})$ par la seconde projection.
  
  \smallskip
  
    Noter que $u$ appartient à $G_u^{^-}$; il s'identifie à l'élément $-1$ de $\GL(V_u^{^-})$.
    
    \smallskip
          
        \noi  On a les inclusions:
          
          \smallskip
  \hspace{8mm}      $G_u^{^+}  \subset    \widetilde{G}_u^{^+}, \ \ G_u^{^-}  \subset    \widetilde{G}_u^{^-} $ \ \ et  \ \ $ G_u^{^+} \times G_u^{^-}   \ \subset \  G_u \ \subset \   \widetilde{G}_u^{^+} \times \widetilde{G}_u^{^-}.$
          
          \medskip
                   
              \noi {\bf Proposition 2.1} (a) {\it $G_u^{^+}$ et $G_u^{^-}$ sont des sous-groupes paraboliques de} $G$.
              
              (b) {\it $G^1_u = G_u^{^+} \times G_u^{^-} $.}
                    
           (c)   {\it $G_u^{^-}$ est engendré par les cubes de $G$ d'extrémité $u$.}
              
              \smallskip
              \noi {\small [Rappelons que $G^1_u$ est le sous-groupe de $G_u$ engendré par les réflexions de $G$ qui commutent à $u$, cf. $\S$1. Un {\it cube} $C$ de $G$ est un sous-groupe abélien engendré par des réflexions; on appelle {\it extrémité} de $C$ l'unique élément de $C$ de degré maximum, cf. [Se 22], 4.1.]}
              
              \smallskip
              
              \noi {\it Démonstration.}
           
    Les groupes $G_u^{^+}$ et $G_u^{^-}$ sont des fixateurs de parties de $V$; cela entraîne que ce sont des {\it sous-groupes paraboliques}; d'où (a). En particulier, ils sont engendrés par des réflexions. D'après la définition de $G^1_u$, on a donc $G^1_u \ \supset G_u^{^+} \times G_u^{^-} $.     
    
    D'autre part, si
   $s \in G_u$  est une réflexion de $G$, on a $s \in G_u^{^+}$ si $\deg(us)=\deg(u)+1$
   et  $s \in G_u^{^-}$ si $\deg(us)=\deg(u)-1$; le groupe $G_u^{^+}$ est donc engendré par les réflexions du premier type, et $G_u^{^-}$ par celles du second type.
      Toute réflexion de $G_u$ est donc contenue dans $G_u^{^+} \times G_u^{^-} $; comme $G^1_u$ est engendré par de telles réflexions, cela démontre l'inclusion $G^1_u \ \subset G_u^{^+} \times G_u^{^-} $.   D'où (b).
      
 Si $C$ est un cube de $G$ d'extrémité $u$, les réflexions appartenant à $C$ sont du second type, donc appartiennent à $G_u^{^-}$, d'où $C \subset G_u^{^-}$. Inversement, toute réflexion de $G_u^{^-}$ appartient à un cube maximal de $G_u^{^-}$; un tel cube a pour extrémité $u$, puisque $u$ est l'élément ``$-1$ '' de $G_u^{^-}$. D'où (c).

\smallskip

   \noi {\bf Corollaire 2.2}. {\it $(V_u^{^+} ,G_u^{^+})$ et $(V_u^{^-} ,G_u^{^-})$ sont
   des couples de Coxeter.}
   
   \smallskip
   \noi {\it Démonstration.} D'après (a), $(V,G_u^{^+})$ est un couple de Coxeter.
   Comme $V = V_u^{^+} \oplus V_u^{^-}$ et que $G_u^{^+}$ opère trivialement sur
   $V_u^{^-}$, il en est de même du couple $(V_u^{^+},G_u^{^+})$. Le cas du couple
   $(V_u^{^-} ,G_u^{^-})$ se traite de manière analogue.
   
           \medskip
        
   \noi {\bf Proposition 2.3}. ([FV 05], prop.7 et [DPR 13], prop.2.2) {\it Le normalisateur de $G_u^{^-}$ dans $G$ est égal à $G_u$.} 

\smallskip

\noi {\it Démonstration.} Il est clair que $G_u^{^-}$ est normal dans $G_u$. Inversement,
soit $g$ un élément de $G$ normalisant $G_u^{^-}$. Comme $u$ est l'unique involution de $G_u^{^-}$ de degré $\deg(u)$, elle est fixée par l'automorphisme intérieur défini par $g$, d'où $g \in G_u$. 

\smallskip

\noi {\it Remarque.} L'énoncé analogue avec $G_u^{^-}$ remplacé par $G_u^{^+}$
n'est pas toujours vrai; il se peut même que $G_u^{^+} = 1$ et $G_u \neq G$; c'est le cas si $G$ est de type $A_2$ et $u$ est une réflexion.

\smallskip
\noi {\bf Proposition 2.4}. {\it La suite exacte $1 \to G_u^1 \to G_u \to \Gamma_u \to 1$ est scindée.}

(Autrement dit, il existe un sous-groupe $X_u$ de $G_u$ tel que $G_u = G^1_u\!\cdot\!X_u$ et $G^1_u \cap X_u=1$.)

\smallskip

\noi {\it Démonstration}. Cela résulte du lemme 2 de [Ho 80],
 appliqué au groupe de Coxeter  $G^1_u$. De plus, la démonstration de [Ho 80] donne une méthode
 pour construire un groupe $X_u$: on choisit une chambre $\cc$ de $G_u^1$ dans $V$,
et on lui associe le sous-groupe $H_{\cc}$ de $\GL(V)$ formé des éléments qui normalisent $G^1$ et qui stabilisent $\cc$. Le normalisateur
de $G_u^1$ dans $\GL(V)$ est le produit semi-direct $G_u^1.H_{\cc}$. On prend alors 
 $X_u =G^1_u \cap H_{\cc}$.

\bigskip

\noi Passons maintenant aux groupes $\widetilde{G}_u^{^+}$ et $\widetilde{G}_u^{^-}$ définis plus haut:

\smallskip
\noi {\bf Proposition 2.5}. {\it On a des suites exactes}:

\smallskip

(a) $1 \ \to \ G_u^{^-} \ \to \ G_u \ \to \widetilde{G}_u^{^+} \ \to \ 1 \ \ \ et \ \  \ 1 \ \to \ G_u^{^+} \ \to \ G_u \ \to  \ \widetilde{G}_u^{^-} \ \to \ 1$.

\smallskip

(b) $1 \ \to \ G_u^{^-} \ \to \ \widetilde{G}_u^{^-} \ \to \ \Gamma_u  \ \to 1  \ \ \ et \ \ \ 1 \  \to \ G_u^{^+} \ \to \ \widetilde{G}_u^{^+} \ \to \ \ \Gamma_u \ \to \ 1$.

\smallskip
\noi{\it Les suites} (b) {\it sont scindées.}

\smallskip

\noi {\it Démonstration}. Le noyau de $G_u \to \GL(V_u^{^+})$ est $G_u^{^-}$; cela entraîne la première suite exacte de (a); la seconde se prouve de la même manière.

D'après (a) on peut identifier $\widetilde{G}_u^{^-}$ à $G_u/G_u^{^+}$. L'homomorphisme $G_u \to \Gamma_u$ est trivial sur $G_u^{^+}$. Il définit donc un homomorphisme $\widetilde{G}_u^{^-} \ \to \ \Gamma_u$ qui est surjectif, et dont le noyau est $G^1_u/ G_u^{^+} = G_u^{^-}$, cf. prop. 2.1 (b). Cela donne la première des suites exactes (b); la seconde se prouve de manière analogue. 

Le fait que ces suites soient scindées résulte du fait analogue pour la suite exacte $1 \to \ G_u^{^+}\times G_u^{^-} \to G_u \to \Gamma_u \to 1,$ cf. prop.2.4.

\smallskip
\noi {\it Remarque}. Les constructions ci-dessus sont des cas particuliers de celles
du {\it lemme de Goursat} ([Se 16], 1.4) qui décrit la structure d'un sous-groupe d'un produit de deux groupes dont les projections sur les deux facteurs sont surjectives. Ici les deux groupes sont
$\widetilde{G}_u^{^+}$ et $\widetilde{G}_u^{^-}$; le sous-groupe est $G_u$.

\smallskip

\noi {\bf Corollaire 2.6}. {\it Dans la suite d'inclusions \
  $G_u^{^+} \times G_u^{^-}   \subset G_u  \subset   \widetilde{G}_u^{^+} \times \widetilde{G}_u^{^-},$
chaque groupe est d'indice $|\Gamma_u|$ dans le suivant.}

\smallskip

  Cela résulte des suites exactes (b). 

\smallskip

\noi {\it Remarque.} Le plus petit des trois groupes du cor.2.6 est normal dans les deux autres. Par contre le groupe du milieu $G_u$ est normal dans le grand seulement si $\Gamma_u$ est abélien. En effet, après passage au quotient par le petit groupe, on obtient
l'inclusion diagonale $1 \ \subset \ \Gamma_u \ \subset \ \Gamma_u \times \Gamma_u$. Or
 la diagonale n'est un sous-groupe normal que si le groupe est abélien.
 
 \medskip
 
  Revenons à la première suite exacte (b): $1  \to  G_u^{^-}  \to  \widetilde{G}_u^{^-}  \to  \Gamma_u   \to 1$. L'action de $\widetilde{G}_u^{^-}$ sur $G_u^{^-}$ par conjugaison donne un homomorphisme $\widetilde{G}_u^{^-} \to \Aut(G_u^{^-})$.
  L'image de cet homomorphisme est contenue dans le sous-groupe $\Aut_c(G_u^{^-})$ de $\Aut(G_u^{^-})$ formé des automorphismes qui transforment
  réflexions en réflexions. Notons $\Out_c(G_u^{^-})$ le quotient de $\Aut_c(G_u^{^-})$ par le sous-groupe des automorphismes intérieurs\footnote{Autre interprétation de $\Out_c(G_u^{^-})$: c'est le groupe des automorphismes du graphe de Coxeter de $G_u^{^-}.$}. Par passage au quotient, on obtient un homomorphisme de $\Gamma_u = \widetilde{G}_u^{^-}/G_u^{^-}$ dans $\Out_c(G_u^{^-})$. 
  
    \medskip
  
  \noi {\bf Proposition 2.7}. {\it L'homomorphisme $\Gamma_u \to \Out_c(G_u^{^-})$ est injectif.}
  
  \smallskip
  
\noi {\it Démonstration.} On utilisera le lemme suivant :

\smallskip
\noi {\bf Lemme 2.8}. {\it Soit $(E,W)$ un couple de Coxeter tel que ${-1 \in W}.$ Tout élément d'ordre fini de $\GL(E)$ qui centralise $W$ appartient
à $W$.}

\smallskip

\noi {\it Démonstration du lemme}. Soit $-1 = \prod_i s_i$ une décomposition de $-1$ en produit de réflexions de $W$ deux à deux distinctes et commutant entre elles. Soit $D_i$ la droite de $E$ sur laquelle $s_i$ opère par $-1$.
On a $E = \oplus_i D_i$. Soit $g \in \GL(E)$ d'ordre fini et centralisant $W$.
Puisque $g$ commute aux $s_i$, il stabilise les $D_i$. Sa restriction à chaque $D_i$ est une homothétie $x \mapsto \varepsilon_i x$, avec $\varepsilon_i \in\R^\times$ d'ordre fini, donc égal à $1$ ou $-1$. Il en résulte
que $g$ est égal au produit des $s_i$ tels que $\varepsilon_i=-1$. En particulier, on a $g\in W$.

\smallskip
\noi{\it Fin de la démonstration de la prop.2.7.} Soit $\gamma$ un élément de $\Gamma_u \to \Out_c(G_u^{^-})$, et soit $g$ un représentant de $\gamma$ dans $\widetilde{G}_u^{^-}$. Supposons que l'image de $\gamma$ dans
$\Out_c(G_u^{^-})$ soit 1. Cela signifie qu'il existe $z \in G_u^{^-}$ tel que
$gxg^{-1} = zxz^{-1}$ pour tout $x \in G_u^{^-}$. L'élément $z^{-1}g$ centralise
  $G_u^{^-}$. D'après le lemme 2.8, appliqué au couple $(V_u^{^-},G_u^{^-})$, on a $z^{-1}g \in G_u^{^-}$, d'où $g\in G_u^{^-}$, i.e. $\gamma = 1$.
  
  \smallskip
  
  \noi {\it Remarque}. On définit de façon analogue un homomorphisme $\Gamma_u \to \Out_c(G_u^{^+})$. Cet homomorphisme est injectif si $-1 \in G$: cela résulte de la prop.2.7, appliquée à $-u$; si $-1 \notin G$, il peut ne pas être injectif.

\bigskip
\noi {\bf Théorème 2.9}. {\it Si le théorème 1.1 est vrai pour $(V,G)$, les couples $(V_u^{^+},\widetilde{G}_u^{^+})$ et
$(V_u^{^-},\widetilde{G}_u^{^-})$ sont des couples de Coxeter. 

En particulier, \ $\widetilde{G}_u^{^+}$ et \ $\widetilde{G}_u^{^-}$ sont des groupes de Coxeter.}

\smallskip
\noi {\it Démonstration}. 

\smallskip
Faisons la démonstration pour  $\widetilde{G}_u^{^+}$; le cas de $\widetilde{G}_u^{^-}$ est analogue.
  Soit $H$ le sous-groupe de $\widetilde{G}_u^{^+}$ engendré par les $V_u^{^+}$-réflexions. Nous devons montrer que
$H = \widetilde{G}_u^{^+}$. Comme ce dernier groupe est un quotient
  de $G_u$, le théorème 1.1 dit qu'il est engendré par les images dans $\GL(V_u^{^+})$ des involutions de degré 1 ou 2 dans $G$. Si $g$ est une involution de $G_u$,
  notons $g^+$ son image dans $\GL(V_u^{^+})$, et $g^-$ son image dans $\GL(V_u^{^-})$. On a $\deg(g)=\deg(g^+)+\deg(g^-)$. Si $\deg(g) \leqslant 2$,
  on a, soit $\deg(g^+) \leqslant 1$, soit $\deg(g^-) = 0$. Dans le premier cas,
  $g^+$ est, soit 1, soit une réflexion dans $\GL(V_u^{^+})$, donc appartient à $H$.
  Dans le second cas, on a $g^-=1$, i.e. $g$ fixe $V_u^{^-}$, donc $g^+$ appartient à
  $G_u^{^+}$, qui est contenu dans $H$, on l'a vu. Cela prouve que $H =\widetilde{G}_u^{^+}$.
  
    \smallskip
  
  \noi {\it Remarque}. Le théorème 2.9 était essentiellement connu, mais dans une formulation différente. On peut le déduire d'un théorème de R.B. Howlett ([Ho 80]) sur les normalisateurs de sous-groupes paraboliques, théorème qui est applicable à $G_u^1$ d'après la prop.2.3. Je dois cette remarque à G. R\"{o}hrle; c'est également lui qui m'a indiqué la prop.2.3.\\
 La démonstration de [Ho 80], comme celle donnée ici, est une vérification cas par cas. Il serait intéressant d'avoir une démonstration directe.

   \bigskip
   \noi  {\bf §3. Détermination des groupes $G^{^+}_u$.}

\smallskip

  Les groupes $G^{^+}_u$ s'obtiennent par une récurrence sur $\deg(u)$ qui permet de passer d'un groupe de Coxeter à un autre de rang inférieur. On se ramène ainsi au cas où $u$ est une réflexion.
  
  \smallskip
   De façon plus précise :

\medskip
\noi {\bf Réduction au cas où $u$ est une réflexion.}
 
\smallskip
\noi {\bf Proposition 3.1}. {\it Soient $v,w$ deux involutions de $G$, commutant entre elles et telles que $\deg(vw)=\deg(v)+\deg(w)$. Alors $G^{^+}_{vw} = (G^{^+}_v)_w^{^+} = (G^{^+}_w)_v^{^+}$}.
 
\smallskip

\noi {\it Démonstration.}  Les hypothèses faites sur $v,w$ équivalent à $V^{^-}_{vw} = V^{^-}_v  \oplus V^{^-}_w.$ Un élément de $G$ appartient à $G^{^+}_{vw}$ si et seulement si il fixe $V^{^-}_v $ et $V^{^-}_w.$ D'où la proposition.

\medskip

\noi {\bf Corollaire 3.2}. {\it Soit $u$ une involution de degré $d$, et soient
$s_1,...,s_d$ des réflexions, commutant deux à deux, telles que $u=s_1\cdots s_d$.
Soit $G(i) \ (i=0,...,d)$ la suite de sous-groupes de $G$ définie par $G(0)=G$ et
$G(i) = G(i-1)^{^+}_{s_i}$. On a $G^{^+}_{u} = G(d)$.}

\smallskip

\noi {\it Démonstration}. Cela résulte de la prop. 3.1 en raisonnant par récurrence sur $d$.

\medskip

\noi {\bf Le cas où  $u$  est une réflexion et où $G$ est cristallographique.}

[Rappelons (cf.[Bo 68], VI.2.5) que $G$ est dit {\it cristallographique} s'il stabilise un réseau de $V$; cela équivaut à dire que $G$ est le groupe de Weyl d'un système de racines de $V$.]

\smallskip

  Supposons que $G$ .soit cristallographique et irréductible. Soit $R$   un système de racines
 de $V$ dont $G$ est le groupe de Weyl, soit $S= \{\alpha_1,...,\alpha_n\}$ une base de $R$ et soit $X$ le graphe de Dynkin  correspondant (celui dont l'ensemble des sommets est $S$).

  Soit $\alpha_0 = - \tilde{\alpha}$ l'opposée de la plus grande racine de $R$
  et soit $X_0 = X \cup \{\alpha_0\}$ le graphe de Dynkin complété (cf. [Bo 68], VI.4.3).
Soit $Y$ le sous-graphe de $X$ obtenu en supprimant les sommets $\{\alpha_i\}$ de $X$ liés à $\alpha_0$ dans $X_0$. Alors :

\smallskip

\noi  {\bf Proposition 3.3}. {\it Soit $s_0$ la réflexion associée à $\alpha_0$. Le
groupe $G^{^+}_{s_0}$ est égal au sous-groupe parabolique $G_Y$ de $G$ de base $Y$.}

\smallskip
\noi {\it Démonstration}. Le groupe $G^{^+}_{s_0}$ est engendré par les réflexions $s_\alpha$ correspondant aux racines positives orthogonales à $\alpha_0$ (pour un produit scalaire défini positif et $G$-invariant, noté $x\!\cdot\! y$). Si l'on écrit
$\alpha$ comme $\sum m_i\alpha_i$, on a $\alpha \!\cdot\! \alpha_0= \sum
\!\cdot\!m_i\alpha_i \!\cdot\! \alpha_0$. Les $m_i$ sont $\geqslant 0$ et les $\alpha_i \!\cdot\! \alpha_0$ sont $\leqslant 0$ ([Bo 68], VI.1.8, prop.8). On a donc $\alpha \!\cdot\! \alpha_0 = 0$ si et seulement si $m_i = 0$ pour tout $i$ tel que $\alpha_i \!\cdot\! \alpha_0 \neq 0$, autrement dit pour tout $i$ tel que $\alpha_i \notin Y$; cela revient à dire que $\alpha$ est une réflexion de $G_Y$, cf. [Bo 68], VI.1.7, cor.4 à la prop.7.
D'où la proposition.

\smallskip

\noi  {\bf Corollaire 3.4}. {\it Supposons que $G$ soit de type impair. Soit $s$ une réflexion de~$G$. Le groupe $G^{^+}_s$ 
est un conjugué du groupe $G_Y$ de la prop.3.3.}

\noi {\small [Rappelons, cf. [Se 22] 1.13, que $G$ est dit de type impair si tous les produits de deux
réflexions sont, soit d'ordre 2, soit d'ordre impair. C'est le cas si $G$ est de l'un des types A, D, E.]}

\smallskip
\noi {\it Démonstration.} Les réflexions d'un groupe de type impair sont conjuguées
entre elles. Donc $s$ est conjuguée de la réflexion $s_0$ de la prop.3. D'où le corollaire.

\smallskip

\noi {\it Exemples}.
  
  \smallskip   Voici trois exemples, qui seront utilisés dans les §§10, 11, 12; les notations sont celles des Tables de [Bo 68], VI.
   
   \smallskip
   
\noi   $(a)$  {\it Type} $E_6$. 
   
  $(a_1)$ Le cas $\deg(u)=1$. Dans le graphe de Dynkin étendu, le sous-diagramme
   $Y$ de la prop.3.3 a pour sommets $\alpha_1, \alpha_3, \alpha_4, \alpha_5, \alpha_6$. Il est de type $A_5$. On a donc $G^{^+}_u \simeq A_5$.
   
   $(a_2)$ Le cas $\deg(u) = 2$. Ecrivons $u$ comme produit de deux réflexions
   $s_1$ et $s_2$, commutant entre elles. D'après $(a_1)$, le groupe $H=G^{^+}_{s_1}$ est de type $A_5$; il contient $s_2$. D'après le cor.3.2, on a $G^{^+}_u = H^{^+}_{s_2}$.
   Comme toutes les réflexions de $H$ sont conjuguées, on en déduit que toutes les
   involutions de degré 2 de $G$ sont conjuguées.
   En appliquant à $H$ le cor.3.4  (ou en raisonnant directement), on voit que $H^{^+}_{s_2}$ est de type $A_3$, et il en est donc de même de $G^{^+}_u$.
   
   $(a_3)$ Le cas $\deg(u) = 3$. Un argument analogue donne à la fois le fait que toutes
   les involutions de degré 3 sont conjuguées et que le groupe $G^{^+}_u$ est de type $A_1$.
   
   $(a_4)$ Le cas $\deg(u)=4$. Même argument : les involutions de degré $4$ sont conjuguées, et le groupe $G^{^+}_u$ est 1.
   
   On peut résumer ce qui précède par une chaîne : $E_6 \longrightarrow A_5 \longrightarrow A_3 \longrightarrow A_1 \longrightarrow1.$
      
      \medskip
    \noi  (b) {\it Type} $E_7$.
      
      La même méthode donne la chaîne $E_7 \longrightarrow D_6 \longrightarrow A_1 \times D_4$, et montre qu'il y a une seule classe de conjugaison d'involutions de degré $2$. Comme un groupe de type $A_1 \times D_4$ a deux types de réflexions,  cette chaîne a deux prolongements possibles, l'un par $D_4$, l'autre
      par $(A_1)^4$; ils correspondent aux deux classes d'involutions de $G$ de degré 3.
      
      \medskip
      
   \noi   (c) {\it Type} $E_8$
      
      Le début de la chaîne est $E_8 \longrightarrow E_7$;  d'après (b), elle se prolonge  par $D_6$, puis par $A_1 \times D_4$, et puis, soit par $D_4$, soit par $(A_1)^4$.

\bigskip
    \noi {\bf 4. Les types $A_1, I_2(m), H_3$ et $H_4$}.
 
 \smallskip
On suppose que le type de $G$ est $A_1, I_2(m), H_3$ ou $H_4$. Soit $u$ une involution de $G$. On se propose de démontrer les théorèmes 1.1 et 1.2 pour le couple $(G,u)$, et de déterminer les groupes $G_u, \Gamma_u,G_u^{^+},..., \widetilde{G}_u^{^-}$ correspondants.

      \medskip
          \noi {\bf Type} $A_1$.
       
       \smallskip
       
       Ici, $G$ est d'ordre 2; l'involution $u$ est, soit $1$, soit $-1$. On a $G_u=G$; si $u=1$, on a $G_u^{^+} = \widetilde{G}_u^{^+}=G$ et $G_u^{^-}=  \widetilde{G}_u^{^-}= 1$; si $u=-1$,
       on a $G_u^{^+} = \widetilde{G}_u^{^+}=1$ et $G_u^{^-}=  \widetilde{G}_u^{^-}= G$.
       Dans les deux cas  $\Gamma_u=1$. Nous résumons
       ceci dans le tableau ci-dessous:

       $$\begin{array}{|c|c|c|c|c|c|c|l}
\cline{1-7}
\hbox{$\deg(u)$} &  |G_u|  &  G_u^{^-} & \widetilde{G}_u^{^-}  &  G_u^{^+} & 
\widetilde{G}_u^{^+} &\gamma_u&\\
\cline{1-7}
0& 2&1& 1 &A_1& A_1& 1 &\\
\cline{1-7}
1& 2&A_1&A_1&1&1&1\\
\cline{1-7}
\end{array}$$
 
   \medskip
       \noi {\bf Type $I_2(m), m$ impair}.
       
       \smallskip
       Le groupe $G$ est diédral d'ordre $2m, m$ impair. Toute involution $u \neq 1$ est une réflexion et son centralisateur est $\{1,u\}$. D'où le tableau :

$$\begin{array}{|c|c|c|c|c|c|c|l}
\cline{1-7}
\hbox{$\deg(u)$} &  |G_u|  &  G_u^{^-} & \widetilde{G}_u^{^-}  &  G_u^{^+} & \widetilde{G}_u^{^+} &\gamma_u & \\ 
\cline{1-7}
0&2m& 1&1&I_2(m)&I_2(m)&1&\\
\cline{1-7}
1& 2&A_1&A_1&1&1&1&\\
\cline{1-7}
\end{array}$$

   \medskip
       \noi {\bf Type $I_2(m), m$ pair}.
       
       \smallskip
       
       Le groupe $G$ est diédral d'ordre divisible par 4. Il contient~-1. Ses
       réflexions forment deux classes de conjugaison, permutées par un automorphisme extérieur; le centralisateur d'une réflexion est le groupe de type $(2,2)$ engendré par cette involution et l'élément $-1$. On en déduit le cas $\deg(u)=1$ du tableau ci-dessous. Le cas où $\deg(u)= 0$ (resp. $2$) est immédiat, puisqu'alors $u=1$ (resp. $-1$).

$$\begin{array}{|c|c|c|c|c|c|c|l}
\cline{1-7}
\hbox{$\deg(u)$} &  |G_u| &  G_u^{^-} & \widetilde{G}_u^{^-}  &  G_u^{^+} & \widetilde{G}_u^{^+} &\gamma_u & \\ 
\cline{1-7}
0&2m& 1&1&I_2(m)&I_2(m)&1&\\
\cline{1-7}
1& 4&A_1&A_1&A_1&A_1&1&\\
\cline{1-7}
2&2m&I_2(m)&I_2(m)&1&1&1&\\
\cline{1-7}
\end{array}$$

   \medskip

\noi {\bf Type $H_3$}.

\smallskip

  Ici, $G = \Alt_5 \times \{1,-1\}$. C'est un groupe de rang 3, contenant $-1$; il y a une seule classe
  d'involutions pour chaque degré $\leqslant 3$. On a le tableau suivant:
$$\begin{array}{|c|c|c|c|c|c|c|l}
\cline{1-7}
\hbox{$\deg(u)$} &  |G_u| &  G_u^{^-} & \widetilde{G}_u^{^-}  &  G_u^{^+} & \widetilde{G}_u^{^+} &\gamma_u & \\ 
\cline{1-7}
0&2^33.5& 1&1&H_3&H_3&1\\
\cline{1-7}
1& 2^3&A_1&A_1&(A_1)^2&(A_1)^2&1&\\
\cline{1-7}
2&2^3&(A_1)^2&(A_1)^2&A_1&A_1&1\\
\cline{1-7}
3&2^33.5&H_3&H_3&1&1&1&\\
\cline{1-7}
\end{array}$$
Les lignes correspondant à $\deg(u) = 0$ ou $3$ sont évidentes. Lorsque $\deg(u)=1$, le groupe $G_u$ est d'ordre $2^3$. On a $G_u^{^-} = A_1$
et $G_u^{^+} = A_1\times A_1$; comme le produit de leurs ordres est égal à celui de
$G_u$, cela montre que $G_u^1 = G_u$ d'où $\Gamma_u=1$. Le cas $\deg(u)=2$ se ramène au précédent en remplaçant $u$ par $-u$, ce qui permute les signes ``$ + $ ''  et ``$- $''.

   \medskip

\noi {\bf Type $H_4$}.

\smallskip
  C'est un groupe de rang 4, contenant $-1$, d'ordre $2^63^25^2$. On l'obtient par ``dédoublement'' à partir de $\Alt_5$, cf. [Se 22], 5.10 et 6.12. Cette construction montre que, pour $d = 0, 1, 2, 3,4$, le nombre des involutions de degré $d$ est respectivement $1, 60, 450, 60,1$, et ces involutions forment une seule classe de conjugaison. On a le tableau suivant :
   
$$\begin{array}{|c|c|c|c|c|c|c|l}
\cline{1-7}
\hbox{$\deg(u)$} &  |G_u|  &  G_u^{^-} & \widetilde{G}_u^{^-} &  G_u^{^+} & \widetilde{G}_u^{^+} &\gamma_u& \\ 
\cline{1-7}
0&2^63^25^2& 1&1&H_4&H_4&1&\\
\cline{1-7}
1& 2^43.5&A_1&A_1&H_3&H_3&1&\\
\cline{1-7}
2&2^5&(A_1)^2&B_2&(A_1)^2&B_2&2\\
\cline{1-7}
3&2^43.5&H_3&H_3&A_1&A_1&1&\\
\cline{1-7}
4&2^63^25^2&H_4&H_4&1&1&1&\\
\cline{1-7}
\end{array}$$

\smallskip

Les cas $\deg(u) = 0$ et $\deg(u)=4$ sont évidents. Le cas $\deg(u)=1$ résulte de ce que le centralisateur d'une réflexion est de type $H_3$; en remplaçant $u$ par $-u$, cela donne le cas $\deg(u)=3$. 

Lorsque $\deg(u)=2$, les groupes $G_u^{^+}$ et $G_u^{^-}$ sont de type $A_1\times A_1$, car sinon ce seraient des groupes diédraux d'ordre $2m$, avec $m$ pair $\geqslant4$, contrairement au fait que $H_4$ est un groupe de type impair, au sens de [Se 22], 1.13 (variante: utiliser le cor.3.4 pour se ramener au type $H_3$). 

Comme
$|G_u|=2^63^25^2/60=2^5$, et que $G_u/G_u^{^-} \simeq G_u^{^+}$, on a $|G_u^{^-}|=8$,
d'où $|\Gamma_u|=2$. Cela justifie la ligne $\deg(u)=2$ du tableau ci-dessus. L'homomorphisme $G_u \to \Gamma_u \simeq \Sym_2$ est donné par l'action de $\Gamma_u$ sur les deux réflexions de produit $u$. Il reste à montrer qu'il existe une involution $g$ de $G_u$, de degré 2, dont l'image dans $\Gamma_u$ est non triviale. Cela résulte d'un énoncé plus général,
démontré au §8. On peut aussi faire un calcul explicite:

\medskip 

Notons  $a,x,y,z$ des réflexions de $G$ réalisant le diagramme de Coxeter

\smallskip

\hspace{46mm}  {$5$}

 \hspace{3cm} $H_4$ : \quad $\circ$------$\circ$------$\circ$------$\circ$ .

\hspace{43mm}$a$\hspace{6mm}$x$\hspace{7mm}$y$\hspace{7mm}$z$

\smallskip

\noi On a $xz=zx, xyx=yxy,yzy=zyz$. Soient $u=xz$ et $g=yuy$. Ce sont des involutions de degré 2.
On a $gxg=z$; en effet, $gxg =yu.yxy.uy=yu.xyx.uy=yzyzy =y.yzy.y=z.$ Ainsi, la conjugaison par $g$ échange
 $x$ et $z$, donc fixe $u$. On a $g\in G_u$, et l'image de $g$ dans $\Gamma_u$ est non triviale.
   
   \bigskip

\noi {\bf 5.  Type $A_{n-1}$}.

  \smallskip

\smallskip
  Dans le cas du type $A$, il est plus commode de décrire $A_{n-1}$ que $A_n$.  Soit $X$ un ensemble fini à $n$ éléments et soit $V_X$ un $\R$-espace vectoriel de base $X$. Le groupe $G=\Sym_X$ des permutations de $X$ opère de façon fidèle sur $V_X$;
  on obtient ainsi un couple de Coxeter $(V_X,G)$ de type $A_{n-1}$; les réflexions  sont les transpositions de $X$. 
  
    \smallskip
  
  L'espace ``$V$'' standard associé à $G$ est l'hyperplan de $V_X$ engendré par les $x-x'$ avec $x,x' \in X$ .
  
  \medskip
  Soit $u\in G$ une involution, autrement dit une permutation de $X$ de carré 1.
Soit $d$ son degré. Soit $Z=X^u$ l'ensemble des points fixes de $u$; le groupe $\{1,u\}$
opère librement sur  $X \sm Z$. Décomposons $X\sm Z$ en deux parties disjointes
$Y, Y'$ telles que $Y'=uY$. On a $d = |Y| = |Y'|$ et $n =a+2d$, où $a=|Z|$. 

Tout élément $g$ de  $G_u$ respecte la décomposition de $X$ en deux parties: $Z$
et $Y \cup Y'$, donc définit une permutation $g_1$ de $Z$ et une permutation $g_2$
de $Y \cup Y'$ commutant à $u$. Inversement, si l'on se donne $g_1,g_2$ vérifiant ces conditions, il lui correspond un élément de $G_u$. Le groupe formé par les
$g_1$ est $\Sym_Z$; il est de type $A_{a-1}$. Celui formé par les $g_2$ est de type
$B_d$, cf. §6. On a donc :

  \smallskip
\noi {\bf Proposition 5.1}. {\it $G_u$ est isomorphe à un groupe de Coxeter de type $A_{a-1} \times B_d$.}

\smallskip

  Soient $x,x'$ deux éléments distincts de $X$. La transposition $\tr_{x,x'}$ appartient
  à $G_u$ si et seulement si l'on a, soit $x,x' \in Z$, soit $x,x' \in Y$ et $x'=ux$.
  Le groupe engendré par les $x,x'$ du premier type est $\Sym_Z$; celui engendré par les $x,x'$ du second type est le produit de $d$ groupes à 2 éléments. Comme les
  transpositions en question engendrent $G^1_u$, on en déduit :
  
  \smallskip
  
  \noi {\bf Proposition 5.2}. {\it Le groupe $G^1_u$ est de type $A_{a-1} \times (A_1)^d$.}
  
    \smallskip
  
  Comme $B_d/(A_1)^d \simeq \Sym_d$, cela entraîne:
  
      \smallskip
\noi{\bf Corollaire 5.3}. {\it On a \ $\Gamma_u \simeq \Sym_d$}.

  \smallskip

Il reste à expliciter les groupes $G_u^{^-}, ...,\widetilde{G}_u^{^+}$. Si $y \in Y$, posons
$y^{+}=y+uy$ et $y^{-}=y-uy$; soit $Y^{+}$ (resp. $Y^{-}$) l'ensemble des $y^{+}$ (resp. des $y^{-}$). Alors  $V_u^{^+}$ a pour base $Z \cup Y^{+}$ et $V_u^{^-}$ a pour base $Y^{-}$. De plus :

\smallskip

\noi (i) Les réflexions de $G_u$ qui fixent  $V_u^{^+}$ sont 
 les transpositions du type $\tr_{y,uy}$, avec $y\in Y$; le groupe $G_u^{^-}$ qu'elles engendrent est de type $(A_1)^Y  \simeq (A_1)^d$.
 
 \smallskip
 \noi (ii) Les réflexions de $G_u$ qui fixent  $V_u^{^-}$ sont 
 les transpositions de $Z$; le groupe $G_u^{^+}$ qu'elles engendrent est $\Sym_Z \simeq \Sym_a$, qui est de type $A_{a-1}$.
 
 \smallskip
 
 \noi (i$'$) Les réflexions de $V_u^{^-}$ qui sont les restrictions d'un élément de $G_u$ sont de deux types :
 
celles de (i), qui changent de signe les $y^{-}$;
  
  celles de la forme $s_{y_1,y_2} = \tr_{y_1,y_2}\tr_{uy_1,uy_2}$, avec $y_1,y_2 \in Y$, qui échangent $y_1^{-}$ et $y_2^{-}$.

  Ces réflexions engendrent un sous-groupe $H$ de  $\widetilde{G}_u^{^-}$ qui est de type $B_d$, donc d'ordre $2^dd!$, cf. §6. Comme $|\widetilde{G}_u^{^-}| = |G_u^{^-}|\cdot|\Gamma_u| = 2^dd!$, on a $H= \widetilde{G}_u^{^-}$. Cela montre que 
  $(V_u^{^-},\widetilde{G}_u^{^-})$ est un couple de Coxeter de type $B_d$, et cela
  montre aussi que $\Gamma_u$ est engendré par les images des $s_{y_1,y_2}$, donc par des images d'involutions de degré 2 de $G_u$. Cela achève la démonstration des th. 1.1 et 1.2 pour $G$.
  
    \smallskip   
\noi (ii$'$) Les réflexions de $V_u^{^+}$ qui sont les restrictions d'un élément de $G_u$ sont celles de (ii), et aussi celles de la forme $s_{y_1,y_2} $, cf.
(i$'$), qui échangent $y_1^{+}$ et $y_2^{+}$.
 On en déduit que $\widetilde{G}_u^{^+}$ est isomorphe à $\Sym_Z \times \Sym_{Y^{+}}$, donc de type $A_{a-1}\times A_{d-1}$.
\medskip

On obtient ainsi le tableau :

 $$\begin{array}{|c|c|c|c|c|c|c|l}
 \cline{1-7}
\hbox{$\deg(u)$} &  |G_u|  &  G_u^{^-} & \widetilde{G}_u^{^-}  &  G_u^{^+} & \widetilde{G}_u^{^+} &\gamma_u&\\
\cline{1-7}
d & 2^dd!a! & (A_1)^d&B_d &A_{a-1} & A_{a-1}\times A_{d-1}& d &\\
\cline{1-7}
\end{array}$$
\noi {\small [Rappelons que $G$ est de type $A_{n-1}$ et que $a$ est le nombre de points fixes de~ $u$.]}
  
  \smallskip
\noi {\small{\it Remarque.} Pour certaines valeurs de $d$, on peut avoir $d-1 =-1$ ou
$n-2d=-1$, ce qui introduit des facteurs $A_{-1}$ dans $G_u^{^+}$  et $\widetilde{G}_u^{^+}$; on les interprète en convenant que $A_{-1} = A_0 = 1$, ce qui est naturel puisque le groupe des permutations d'un ensemble à $0$ ou $1$ élément est égal à 1.  Dans les tableaux relatifs aux types $B_n$ et $D_n$, on rencontre aussi $B_0, B_1, D_0, D_1$; on convient que $B_0=1, B_1=A_1, D_0=D_1=1$.}

   \bigskip

\noi {\bf 6.  Type $B_n$}.

\smallskip

  Soit $n$ un entier $> 0$. (On pourrait même supposer $n>2$, car $B_1=A_1$
  et $B_2=I_2(4)$, et ces cas ont été traités au §4.)

  \medskip 
  
\noi {\bf Rappels}.

\smallskip
  
  Soit  $Z$ est un ensemble fini à $2n$ éléments, et soit $\varepsilon$ une permutation de $Z$ de carré $1$ sans point fixe. Soit $Y$ le quotient
  de $Z$ par l'action du groupe $C=\{1,\varepsilon\}$. On a $|Y|=n$. 
  Soit $G=\Sym_{Z,Y}$ le groupe des permutations de $Z$ commutant à $\varepsilon$, autrement dit le groupe d'automorphismes du diagramme $Z \to Y$. On a une suite exacte $$1 \to C^Y \to G \to {\rm Sym} _Y \to 1,$$
 où $C^Y$ est le groupe des applications de $Y$ dans $C$ (i.e. un produit de $n$ copies de $C$ indexées par $Y$). Cette suite est scindée. On a
  $|G| = 2^n n!$ .
  
  Soit $V_Z$ un $\R$-espace vectoriel de base $Z$. Le groupe $G$ opère sur 
  $V_Z$, et stabilise le sous-espace formé des éléments invariants par $\varepsilon$,
  espace qui s'identifie à $V_Y$; soit $V= (1-\epsilon)V_Z$ l'espace formé par les anti-invariants de $\varepsilon$; on a $V_Z = V_Y \oplus V$. L'action de $G$ sur 
  $V$ est fidèle.  {\it Le couple $(V,G)$ est un couple de Coxeter de type $B_n$}.
  
    Il y a deux classes de réflexions : les {\it courtes} qui sont des transpositions de la forme $\tr_{z,\varepsilon z}$, avec  $z \in Z$, et les {\it longues} qui sont de la forme $\tr_{z,z'}\tr_{\varepsilon z, \varepsilon z'}$, avec  $z,z' \in Z$ et $ z' \neq z, \varepsilon z $. Les premières sont des permutations impaires de $Z$, et les secondes sont des permutations paires. Avec les notations de [Bo 68], VI.4.5, ces réflexions correspondent aux racines $±\epsilon_i$ et $±\epsilon_i±\epsilon_j \ (i \neq j)$.
    
    \smallskip
\noi {\it Remarque}.  Le groupe $G$, vu comme sous-groupe de
$\Sym_Z$ n'est pas engendré par des transpositions si $n>1$; mais il est engendré par des transpositions et des produits de deux transpositions, autrement dit par des involutions de $\GL(V_L)$ de degré 1 ou 2.

  \medskip
  
  \noi {\bf Les groupes $G_u, \ G^1_u$ et $ \Gamma_u$ associés à une involution $u$}.
  
  \smallskip
    Soit $u$ une involution de $G$, autrement dit une permutation de $Z$ de carré 1 qui commute
    à $\varepsilon$. Lorsque $u=1$ ou $u=\varepsilon$, on a $G_u=G$. Supposons que
    $u \neq 1, \varepsilon$. Soit $\Delta=\langle u,\varepsilon\rangle$ le groupe d'ordre 4 engendré par $u$ et $\varepsilon$. L'action de $\Delta$ sur $Z$ donne une partition de $Z$ en trois sous-ensembles :
    
    \smallskip
    
  \   $Z_u \ = $ ensemble des $z\in Z$ tels que $uz =\varepsilon z$;
    
   \ $Z'_u \ =  $ ensemble des $z\in Z$ tels que $uz=z$;
    
   \  $Z''_u \ = $   ensemble des $z\in Z$ tels que $uz \neq z, \varepsilon z$.
    
    \smallskip
    \noi Ces ensembles sont stables par $\Delta$, donc par $\varepsilon$; soient $Y_u, Y'_u,Y''_u$ leurs
    images dans $Y$: on obtient ainsi une partition $Y = Y_u \cup Y'_u \cup Y''_u$. 
    L'ensemble des points de $Y$ fixés par $u$ est $Y_u \cup Y'_u$. Le groupe $\{1,u\}$ opère librement sur $Y''_u$; soit $T_u$ le quotient de $Y''_u$ par cette action. Cela donne le diagramme :

 \bigskip

\hspace{25mm}$Z \ = \quad Z_u \quad \cup \quad  Z'_u \quad  \cup \quad  Z''_u$

\hspace{38mm}$ \downarrow \hspace{7mm} \qquad \downarrow \hspace{13mm} \downarrow$

\hspace{25mm}$Y \ = \quad Y_u \quad  \cup  \quad Y'_u \quad  \cup \quad \ Y''_u$

$\hspace{69mm} \downarrow$

$\hspace{68mm} T_u$.

\bigskip
    Posons \ $ a= |Y_u|, \ a'= |Y'_u|, \ b= |T_u| = \frac{1}{2}|Y''_u|.$ On a $n = a+a'+2b$ et $\deg(u)=a+b$. Les entiers $a$ et $b$ caractérisent la classe de conjugaison de l'involution $u$, et peuvent être donnés arbitrairement pourvu que $a+2b \leqslant n$.
      
    \smallskip
    Le quadruplet $(Z,\varepsilon,Y,u)$ est réunion disjointe de trois quadruplets correspondant aux trois composantes de $Y$ que l'on vient de définir.
    Cette décomposition est stable par le groupe $G_u$. Plus précisément, $G_u$
    est produit direct de trois facteurs:

  \smallskip
  (i) Le premier facteur est $\Sym_{Z_u,Y_u}$ est d'ordre $2^aa!$; sa contribution à
  $G_u^{^-}$ et à  $\widetilde{G}_u^{^-}$ est $\Sym_{Z_u,Y_u}$; celle à $G_u^{^+}$ et à  $\widetilde{G}_u^{^+}$ est $1$; celle à $\Gamma_u$ est $1$. 
  
  \smallskip
  (ii) Le second facteur est $\Sym_{Z'_u,Y'_u}$ est d'ordre $2^{a'}a'!$; sa contribution à
  $G_u^{^-}$ et à  $\widetilde{G}_u^{^-}$ est $1$;  celle à $G_u^{^+}$ et à  $\widetilde{G}_u^{^+}$ est $\Sym_{Z'_u,Y'_u}$; celle à $\Gamma_u$ est $1$.

\smallskip
      
    (iii) Le troisième facteur est le groupe $\Sym_\Delta(Z''_u)$ des permutations de $Z''_u$ qui commutent à l'action de $\Delta$. Il est produit semi-direct des deux groupes suivants :
    
  $\bullet$     le groupe $\Delta^{T_u}$ des applications de $T_u$ dans $\Delta$; c'est un sous-groupe normal d'ordre $4^b$;
            
$\bullet$   le groupe des $\Delta$-automorphismes de $Z''_u$ qui stabilisent une partie $S$ rencontrant chaque fibre de $Z''_u \to T_u$ en un point et un seul; il
 est isomorphe à $\Sym_{T_u}$. 
 
 On peut donc écrire le troisième facteur sous la forme $ \Delta^{T_u}\!.\Sym_{T_u}.$
 Son ordre est $4^bb!$ . Si $b>1$, ce n'est pas un groupe de Coxeter.
     
     \smallskip
      
   \noi   En résumé :
      
      \smallskip
       \noi {\bf Proposition 6.1}. {\it On a} \ $G_u = {\rm \Sym}_{Z_u,Y_u} \times {\rm \Sym}_{Z'_u,Y'_u} \times  \Delta^{T_u}\!.\Sym_{T_u}.$
        
       \smallskip
       
       Les deux premiers facteurs de $G_u$ sont engendrés par des réflexions;
       ils sont donc contenus dans $G^1_u$. Il n'en est pas de même du troisième facteur si $b>1$:

\smallskip
      
      \noi {\bf Proposition 6.2}. {\it On a \ $ G^1_u \ \cap \Delta^{T_u}\!.\Sym_{T_u}  = \ \Delta^{T_u}$.}

         \smallskip
         
         \noi {\it Démonstration}. Il suffit de montrer que toute réflexion $s$ de $G_u$
         qui fixe $Z$ et $Z'$ appartient au groupe $\Delta^{T_u}$. Soit $I_s$ l'ensemble des points de $Z''$ qui ne sont pas fixés par $s$; c'est un ensemble à $2$ ou à $4$
         éléments, on l'a vu. Or $I_s$ est stable par $\Delta$, et les orbites de $\Delta$
         dans $Z''$ sont d'ordre 4. Donc $I_s$ est une orbite de $\Delta$, ce qui entraîne
       que $s$ appartient à $\Delta^{T_u}$.
       
       \smallskip
   
\noi {\bf Proposition 6.4}.
(a) {\it L'action de $G_u$ sur $T_u$ définit par passage au quotient un isomorphisme de \ $\Gamma_u$ sur $\Sym_{T_u}$.}

(b) {\it Toute transposition de $\Sym_{T_u}$ est image d'une involution de $G_u$
de degré~$2$.}

\smallskip

\noi {\it Démonstration de} (a). Cela résulte des prop. 6.1 et 6.2 puisque $\Gamma_u=G_u/G^1_u$.

\noi {\it Démonstration de} (b). Soient $t,t' \in T_u$, avec $t \neq t'$, et soient $z,z'$ des représentants de $t,t'$ dans $Z''_u$. Soient $g,h$ les réflexions de $G$ données par $g =\tr_{z,z'}\tr_{\varepsilon z,\varepsilon  z'}$ et $h = \tr_{uz,uz'}\tr_{\varepsilon  uz,\varepsilon  uz'}$. Ces réflexions commutent, et l'on a $ugu=h$. On a $gh \in G_u$ et l'image de $gh$  dans $\Sym_{T_u}$ est la transposition $\tr_{t,t'}$. D'où (b).

\medskip

\noi {\bf Corollaire 6.5}. {\it Les théorèmes 1.1 et 1.2 sont vrais pour $G$.}

\smallskip C'est clair.

  \smallskip
  
On peut résumer les résultats obtenus de la façon suivante:
  
  \smallskip
   
    \noi {\bf Proposition 6.6}. {\it On a} :
    
     \smallskip
     $G_u =\Sym_{Z_u,Y_u} \times  \Sym_{Z'_u,Y'_u} \times \ \Delta^{T_u}\!. \Sym_{T_u}  \ \simeq \ B_a\times B_{a'}\times \Delta^b.\Sym_b$;
     
      \smallskip
     $G^1_u = \Sym_{Z_u,Y_u} \times  \Sym_{Z'_u,Y'_u} \times \hspace{6mm} \Delta^{T_u}\hspace{7mm} \simeq \ \  B_a\times B_{a'}\times \Delta^b$;
     
      \smallskip
     $\Gamma_u = \Sym_{T_u} \simeq \Sym_b$ ; $\gamma_u = |T_u| = b$.
   
    \bigskip

  \noi {\bf Les groupes $G^{^+}_u, \widetilde{G}^{^+}_u, G^{^-}_u, \widetilde{G}^{^-}_u$ associés à $u$}.
  
  \smallskip
  
    La décomposition de $G$ en produit de trois facteurs entraîne une décomposition du même type pour les groupes $G^{^+}_u, ...,\widetilde{G}^{^-}_u$. Nous avons donné plus haut le cas des deux premiers facteurs. Pour le troisième facteur, on a:
    
    \smallskip
    
    \noi {\bf Lemma 6.7}. {\it Les troisièmes facteurs des groupes $G^{^+}_u$  et 
    $G^{^-}_u$ sont de type $(A_1)^b$. Ceux des groupes $\widetilde{G}^{^+}_u$ et
   $ \widetilde{G}^{^-}_u$ sont de type $B_b$.}
        
        \smallskip
        
        \noi {\it Démonstration}. Puisque cet énoncé ne concerne que le troisième facteur,
        on peut supposer que les deux premiers sont triviaux, i.e. que $a=a'=0$
        et $Z=Z''_u$. Dans ce cas, les involutions $u$ et $\varepsilon u$ sont conjuguées, ce qui entraîne que $G^{^+}_u \simeq
    G^{^-}_u$ et $\widetilde{G}^{^+}_u \simeq \widetilde{G}^{^-}_u$. Notons ces groupes $H$ et $\widetilde{H}$. D'après la prop.6.2, on a $H \times H \simeq \Delta^b$, ce qui entraîne que $H$ est un groupe abélien élémentaire d'ordre $2^b$;
    comme c'est un groupe de Coxeter, il est isomorphe à $(A_1)^b$. Un argument analogue montre que $\widetilde{H}\times\widetilde{H} \simeq \Delta^b.\Sym_b$.
    En particulier $\widetilde{H}$ est d'ordre $2^dd!$. Or, il contient $H$ comme sous-groupe normal. Cela entraîne que c'est un groupe de Coxeter de type $B_b$,
    en vertu du lemme suivant:
    
    \smallskip
    
    \noi {\bf Lemme 6.8}. {\it Soit $(E,H)$ un couple de Coxeter. Soit $e = \dim E$.
    Supposons que $H$ soit de type $(A_1)^e$. Soit $H'$ un sous-groupe fini de $\GL(E)$ qui normalise $H$ et qui est d'ordre $2^ee!$ . Alors $(E,H')$ est un couple
    de Coxeter de type $B_e$. } 
    
    \smallskip
    
    \noi {\it Démonstration du Lemme 6.8}. L'action de $H$ décompose $E$ en somme directe de droites $D_1,...,D_e$. Comme $H'$ normalise $H$, il permute les $D_i$.
   Soit   $\langle x\!\cdot\!y\rangle$ un produit scalaire défini positif sur $E$ invariant par $H'$,
   et soit $Z$ l'ensemble des $z \in D_1 \cup ... \cup D_e$ tels que $\langle z\!\cdot\!z\rangle=1$. On a $|Z| = 2e$, et l'application $z \mapsto -z$ est une permutation $\varepsilon$ d'ordre 2 de $Z$ sans point fixe. Le groupe $H'$ stabilise $Z$ et commute à $\varepsilon$. On obtient ainsi un homomorphisme injectif de $H'$
   dans le groupe de Coxeter de type $B_e$ défini par $(Z, \varepsilon)$; comme les deux groupes ont le même ordre, cet homomorphisme est un isomorphisme.              
    
    \medskip
    
   \noi On obtient finalement le tableau  :

   $$\begin{array}{|c|c|c|c|c|c|c|l}
   \cline{1-7}
\hbox{invariants} &  |G_u|  &  G_u^{^-} & \widetilde{G}_u^{^-} &  G_u^{^+} &\widetilde{G}_u^{^+} &\gamma_u & \\ 
\cline{1-7}
a, a', b & 2^na!a'!b!&B_a\times (A_1)^b & B_a\times B_b&B_{a'}\times (A_1)^b&B_{a'}\times B_b&b&\\
\cline{1-7}
\end{array}$$

  \  \smallskip
\noi {\bf 7.  Type $D_n$}.

\smallskip
Conservons les notations $(Z, \varepsilon, Y)$ du §6. Soit $G'=B_{Z,Y}$ et soit $G= D_{Z,Y}$ le sous-groupe d'indice $2$ de $G'$ formé des éléments $g$ qui sont des permutations {\it paires} de $Z$, i.e. $\sgn_Z(g)=1$. Le couple $(V_{Z,Y},G)$ est un couple de Coxeter de type $D_n$. Les réflexions de $G$ sont les réflexions longues de $G'$.

  Soit $u$ une involution de $G$, et soient $a,a',b$ ses invariants au sens du §6. Le fait que $u$
  appartienne à $G$ équivaut à $a \equiv 0$ (mod 2). Deux involutions de mêmes invariants sont conjuguées, sauf dans le cas $a=a'=0$ où il y a deux classes de conjugaison.
  
  Le groupe $G_u$ est le sous-groupe d'indice $\leqslant 2$ de $G'_u$ formé des éléments $x$ tels que $\sgn_Z(x)=1$. La décomposition de $G'_u$ donnée dans la prop. 6.1 est :
  
  \smallskip
  
  $G'_u = {\rm \Sym}_{Z_u,Y_u} \times {\rm \Sym}_{Z'_u,Y'_u} \times \Aut_\Delta(Z''_u).$
    
\smallskip

Les éléments de $\Aut_\Delta(Z''_u)$ sont de signature 1. La condition $\sgn_Z(x)=1$
ne porte donc que sur les deux premières composantes de $x$. D'où :

  \smallskip
  
  \noi {\bf Proposition 7.1.} {\it Soit $H_u$ le sous-groupe de ${\rm \Sym}_{Z_u,Y_u} \times {\rm \Sym}_{Z'_u,Y'_u}$ formé des couples $(g,g')$ tels que $\sgn_{Z_u}(g)=
  \sgn_{Z'_u}(g')$. On a $G_u = H_u \times \Aut_\Delta(Z''_u)$.}
  
 \noi {\small[Rappelons que $\Delta=\langle u,\varepsilon\rangle$.]}
 
 \smallskip
 
 Il y a quatre possibilités pour $(a,a')$ :
 
 \smallskip
 (i)  $a=a'=0$, i.e. $2b=n$. On a alors $H_u=1$ et $G_u = \Aut_\Delta(Z''_u) = G'_u$. L'ordre de $G_u$ est $2^nb!$, on a $\Gamma_u = \Sym_{T_u}$ et $\gamma_u = b$. Les groupes
 $ G_u^{^-}$  et  $G_u^{^+}$  sont isomorphes à $A_1^{T_u}$; les groupes $ \widetilde{G}_u^{^-}$  et  $\widetilde{G}_u^{^+}$ sont de type $B_b$.
  
  \smallskip
  (ii) $a=0, a' >0$. Le premier facteur de $H_u$ est 1;  le second est $D(Z'_u,Y'_u)$,
  qui est de type $D_{a'}$. On a  $G_u= D(Z'_u,Y'_u)\times \Aut_D(Z''_u)$.
  La situation est la même que pour $G'$, avec $B_{a'}$ remplacé par $D_{a'}$.
  On a $\Gamma_u = \Sym_{T_u}$ et $\gamma_u = b$.
  
  \smallskip
  (iii) $a>0$ et $a'=0$: comme dans le type (ii), avec $a$ et $a'$ permutés, ainsi que
  $(Z,Z')$ et $(Y,Y')$. Ici encore  $\Gamma_u = \Sym_{T_u}$ et $\gamma_u=b$.
  
  \smallskip
  
  (iv) $a>0$ et $a'>0$. Soit $H^1_u = D_{Z_u,Y_u} \times D_{Z'_{u},Y'_{u}}$. C'est
  un sous-groupe d'indice 2 de $H_u$ qui est engendré par des réflexions. Inversement, toute réflexion de $H_u$ appartient à $H^1_u$ car c'est
  transposition de $Y_u \cup Y'_u$ qui stabilise à la fois $Y_u$ et $Y'_u$, donc qui est une transposition, soit  de $Y_u$, soit de $Y'_u$. D'autre part, le groupe engendré
  par les réflexions de $\Aut_\Delta(Z''_u)$ est le groupe $\Delta^{T_u}.$ On en conclut que $G^1_u = H^1_u \times \Delta^{T_u}$ et que le groupe $\Gamma_u= G_u/G^1_u)$ est égal au produit de $\Sym_{T_u}$ par $H_u/H^1_u$ qui est d'ordre 2. C'est le cas, mentionné dans le th. 2.2, où $\Gamma_u$ {\it n'est pas un
  groupe symétrique} (sauf si $b=0$ ou $1$).  
  
    \medskip

 La détermination des groupes $G_u^{^-}, G_u^{^+},...$ résulte de celle des groupes correspondants pour le type $B_n$. Plus précisément :
 
   Les composantes dépendant de l'invariant ``$b$ '' sont les mêmes que pour le type $B_n$; dans les autres, certains groupes $B_a$ ou $B_{a'}$ sont remplacés par $D_a$ ou $D_{a'}$ respectivement.
   
On obtient ainsi le tableau :
$$\begin{array}{|c|c|c|c|c|c|c|l}
\cline{1-7}
\hbox{invariants} \ \ \ & |G_u| & G_u^{^-} & \widetilde{G}_u^{^-} &  G_u^{^+} & \widetilde{G}_u^{^+}&\!\gamma_u\!& \\ 
\cline{1-7}
0,0,b\ \hspace{10mm} &\!2^nb!&(A_1)^ b   & B_b&   (A_1)^b& B_b&\!b\!&\\
\cline{1-7}
0,a',b \ \ a'>0 \!& \!2^{n-1}a'!b!\!&\!  \!(A_1)^b\!&\! \!B_b\!&\!D_{a'}\times \!(A_1)^b\!&\!D_{a'}\!\times\!B_b\!&\!b\!&\\
\cline{1-7}
a,0,b \ \ \ a>0\!& \!2^{n-1}a!b!\!&\!D_a\! \times \!(A_1)^b\!&\!D_a\times \!B_b\!&\!(A_1)^b\!&B_b\!&\!b\!&\\
\cline{1-7}
\! a,a',b \ \ aa'>0 \! &\!2^{n-1}a!a'!b!\!&\!D_a\! \times \!(A_1)^b\!&\!B_a\times \!B_b\!&\!D_{a'}\times \!(A_1)^b\!&\!B_{a'}\!\times\!B_b\!&\!b,2\!&\\
\cline{1-7}
\end{array}$$

\bigskip
    
\noi {\bf 8. Résultats auxiliaires sur les groupes $G_u^{^-}$ et $\widetilde{G}_u^{^-}$}.

\medskip

  Ces résultats seront utilisés dans les trois sections suivantes. On  note $d$ le degré de l'involution $u$. On suppose que le type de $G$ est $H_4,E_6, E_7$ ou $ E_8$. On s'intéresse aux deux propriétés suivantes:    
 \medskip
 
 (i) $u$ {\it est l'extrémité d'un seul cube}, i.e. sa décomposition en produit de $d$
 réflexions est unique, à permutation près. Cela équivaut à $G_u^{^-} \simeq (A_1)^d$.
 
 \medskip
 (ii) {\it Il existe un $\cc$-sous-groupe $H$ de $G$ de type $A$ tel  que $u\in H$.} 
 
  \medskip
{\bf Proposition 8.1}. {\it Si les propriétés} (i) {\it et} (ii) {\it sont satisfaites, les théorèmes 1.1 et 1.2 sont vrais pour le couple $(G,u),$ alors
le groupe $\widetilde{G}_u^{^-}$ est de type $B_d$, le groupe $\Gamma_u$ est isomorphe à $\Sym_d$ et il 
 est engendré par les images des involutions de $G_u$ de degré $2$.}
 
 \medskip
\noi {\it Démonstration}. Soit $H$ un sous-groupe de $G$ satisfaisant à (ii). D'après le §5, la prop.8.5 est vraie si  $G=H$. On va se ramener à ce cas. D'après (i), on a $G_u^{^-} =  H_u^{^-}$. Le groupe $\widetilde{G}_u^{^-}$ contient $\widetilde{H}_u^{^-}$, qui est de type
$B_d$. Cela montre que $\Gamma_u$ contient un sous-groupe isomorphe à $\Sym_d$. D'autre part $\Gamma_u$
   est isomorphe à un sous-groupe de $\Out_c(G_u^{^-}) \simeq  \Out_c((A_1)^d) \simeq\Sym_d$, cf. prop.2.7. On a donc
$\widetilde{G}_u^{^-}=\widetilde{H}_u^{^-}$, ce qui démontre la proposition.

\bigskip

{\bf Proposition 8.2}. {\it Les propriétés} (i) {\it et} (ii) {\it sont satisfaites pour $d \leqslant 2$ lorsque 
$G$ est de type $H_4$ et pour $d \leqslant 3$ lorsque $G$ est de type $E_6, E_7$ ou $E_8$. La propriété} (ii) {\it est satisfaite pour $d \leqslant 4$ lorsque $G$ est de type $E_8$.}

\smallskip
\noi {\it Démonstration}. L'hypothèse (i) est satisfaite puisque les seuls groupes
de Coxeter de rang $\leqslant 3$, de type impair, et contenant $-1$, sont des puissances de $A_1$.

 Pour (ii), et $G$ de type $H_4$, on remarque que le diagramme de $G$
contient un sous-diagramme de type $A_3$. Or un groupe de type $A_3$ contient
des involutions de tout degré $\leqslant 2$. Comme les involutions de $G$ de même degré sont conjuguées entre elles, cela entraîne (ii). 

Le même argument s'applique à $G$ de type $E_6$, car son diagramme contient un sous-diagramme de type $A_5$;
il s'applique aussi au type $E_8$, ainsi qu'à $E_7$ si $d\leqslant 2$. 

Dans le cas de $E_7$, pour $d=3$,  il y a deux classes d'involutions, cf. [Se 22], 7.5: celles de type {\it triangle} et celles de type {\it droite}. Pour les traiter, choisissons des réflexions $s_1,...,s_7$ correspondant au diagramme de Coxeter de $E_7$ :

\hspace{5mm}$s_{_1} \hspace{6mm}s_{_3}\hspace{6mm}s_{_4}\hspace{6mm} s_{_5}\hspace{6mm} s_{_6}\hspace{6mm} s{_7}$

\hspace{5mm}$\circ$-------$\circ$------$\circ$-------$\circ$------$\circ$------$\circ$

\hspace{23mm} $\mid$

\hspace{23mm} $\circ \ s_{_2}$

Posons $u= s_3s_5s_7$ et $u'=s_2s_5s_7$; ce sont des involutions de degré $3$. D'après [Se 22],  {\it loc.cit.}, un produit $s_as_bs_c$ est du type triangle si et seulement si il existe $m≠ a,b,c$ tel que $s_m$ soit adjacent à un et un seul des
$s_a,s_b,s_c$\footnote{Dans [Se 22],  cette condition est exprimée en termes des racines $\alpha_a,\alpha_b,\alpha_c$ associées à $s_a,s_b,s_c$ : le produit $s_as_bs_c$ est du type triangle si et seulement si $\frac{1}{2}(\alpha_a+\alpha_b+\alpha_c)$ n'appartient pas au réseau des poids.}. Dans le cas de $u$, l'entier $m=1$ répond à cette condition; dans le cas
de $u'$, aucun $m$ n'est possible; ainsi, 
 $u$ est du type triangle et $u'$ du type droite.  Il suffit donc de vérifier la condition (ii) pour  $u$ et pour $u'$. Pour $u$, on prend le sous-groupe de type $A_5$ engendré par $s_3,s_4,s_5,s_6,s_7$ ; pour $u'$, on prend celui engendré par $s_2,s_4,s_5,s_6,s_7$. 
 
 Le cas de $E_8$ est analogue au précédent. Il y a deux classes d'involutions de degré $4$, celles de type {\it rectangle} et celles de type {\it tétraèdre}. On les distingue
 de la manière suivante : on écrit l'involution $u$ comme produit de quatre
 réflexions commutant deux à deux, et correspondant à des racines $x,y,z,t$.
 Alors $u$ est de type rectangle si $x+y+z+t \in 2R$, où $R$ désigne le réseau des racines; sinon, $u$ est de type tétraèdre. Soient $s_1,...,s_8,s_0$ des réflexions correspondant au diagramme étendu de $E_8$
 ($s_0$ correspondant à la plus grande racine): 
 
\vspace{5mm}
\hspace{5mm}$s_{_1} \hspace{6mm}s_{_3}\hspace{6mm}s_{_4}\hspace{6mm} s_{_5}\hspace{6mm}s_{_6}\hspace{6mm} s_{_7}\hspace{6mm} s_{_8}\hspace{6mm} s_{_0}$

\hspace{5mm}$\circ$-------$\circ$------$\circ$-------$\circ$------$\circ$------$\circ$------$\circ$------$\circ$

\hspace{23mm} $\mid$

\hspace{22mm} \ $\circ \ s_{_2}$

\smallskip
Prenons 
$ u = s_2s_5s_7s_0$ et $u' = s_1s_3s_5s_7$. On vérifie par la même méthode que
 pour $E_7$ que $u$ est du type rectangle et $u'$ du type tétraèdre. Il suffit donc
 de vérifier (i) pour $u$ et pour $u'$: pour $u$ (resp. pour $u'$) on prend le sous-groupe de type $A_7$
 engendré par les $s_i, i \neq 1,3$ (resp. $i\neq0,2$).

\bigskip

\noi {\bf 9.  Type $E_6$}.

\smallskip

C'est un groupe de rang $6$ qui ne contient pas $-1$.
Il y a une seule classe de conjugaison d'involutions pour chaque degré $d\leqslant 4$. Le tableau correspondant est :

\smallskip
$$\begin{array}{|c|c|c|c|c|c|c|l}
\cline{1-7}
\hbox{$\deg(u)$} &  |G_u| &  G_u^{^-} & \widetilde{G}_u^{^-}  &  G_u^{^+} &\widetilde{G}_u^{^+} & \gamma_u& \\ 
\cline{1-7}
0 & 2^73^45 & 1 & 1 \ & \  E_6 & E_6 & 1 &\\
\cline{1-7}
1 & 2^53^25 & A_1 & A_1& A_5& A_5  & 1 & \\
\cline{1-7}
2 & 2^63& (A_1)^2 & B_2 & 
A_3 &   A_1\times A_3  & 2&\\
\cline{1-7}
3 &  2^53 & (A_1)^3  & B_3\  &  A_1  & \ \ A_1 \times A_2 \  \   & 3 &\\
\cline{1-7}
4 & 2^73^2 & D_4 & F_4 & 1 & A_2 & 3 &\\
\cline{1-7}
\end{array}$$

\bigskip
\noi {\it Vérification du tableau.}

\smallskip

Les cas $d=0$ et $d=1$ sont immédiats. 

Pour $d= 2$ ou $3$, les prop.8.1 et 8.2 montrent que le
th.1.1. et le th.1.2 sont vrais pour $(G,u)$, que $G_u^{^-} \simeq (A_1)^d$, $\widetilde{G}_u^{^-}\simeq B_d$  et que $\Gamma_u \simeq \Sym_d$.

Dans le cas $d=2$, l'ordre de $G_u$ est $2^63$ et l'on a vu au §3 que $G_u^{^+}$
est de type $A_3$.
Le groupe  $\widetilde{G}_u^{^+}$ contient $G_u^{^+}$ comme sous-groupe d'indice $2$; de plus c'est un groupe de type cristallographique; la seule possibilité est
qu'il soit de type $A_1 \times A_3$.

Dans le cas $d=3$, un raisonnement analogue montre que l'ordre de $G_u^{^+}$ est $2$,
donc que ce groupe est de type $A_1$; quant à  $\widetilde{G}_u^{^+}$ , c'est un groupe de rang $\leqslant 3$ et d'ordre $2^23$, qui normalise un groupe de type $A_1$; son type est donc $A_1 \times A_2$.

Le cas $d=4$ est traité dans [Se 22], 3.19 (et c'est lui qui est à l'origine du présent travail).

  \bigskip

\noi {\bf 10.  Type $E_7$}.

\smallskip

Il y a une seule classe d'involutions de degré $0, 1, 2, 5,6,7$, deux classes
de degré 3 et deux classes de degré 4; ces dernières sont notées $3, 3', 4, 4'$ dans le tableau ci-dessous.
  Rappelons comment on caractérise les deux classes de degré 3, cf. [Se 22], 7.5. Notons $R$ le réseau des racines de $E7$, $P$ le réseau des poids, et $V_6$ le $\F_2$-espace vectoriel $R/2P$; cet espace est muni d'une forme bilinéaire alternée non dégénérée. Les réflexions de $E_7$ correspondent bijectivement aux éléments non nuls de $V_6$. Soit une involution de degré 3; décomposons  $u$ en produit de trois réflexions $s,s',s''$ commutant entre elles (ce qui est unique, à permutation près); ces réflexions donnent trois éléments non nuls $x,x',x''$ de $R/2P$, deux à deux orthogonaux. Si la somme $x+x'+x''$ est $0$,
 $u$ est du type {\it }droite; si elle ne l'est pas, $u$ est du type {\it triangle}. Il n'est pas difficile de compter combien il y a d'involutions de chaque type : pour le cas des droites, il y a $2^6-1$ possibilités pour $x$, $2^5-2$ possibilités pour $x'$ et une seule possibilité pour $x''$. Comme
 chaque $u$ est obtenu $6$ fois, le nombre des involutions du type droite est
 $(2^6-1)(2^5-2)/6 = 3^25.7$ et $G_u$ est d'ordre
 $2^{10}3^45.7/3^25.7=2^{10}3^2$: c'est le type 3 du tableau ci-dessous. Un calcul analogue montre que le nombre des involutions de type triangle est $(2^6-1)(2^5-2)(2^4-2^2)/6 = 2^23^35.7$ et que $G_u$ est d'ordre $2^83$;  c'est le type $3'$ du tableau.
   Si $u$ est de degré 4 (resp. $4'$), on dit que $u$ est de type 4 si $-u$ est de type 3 
 (resp. $3'$).
 
\smallskip 

$$\begin{array}{|c|c|c|c|c|c|c|l}
\cline{1-7}
\hbox{$\deg(u)$} &  |G_u|  &  G_u^{^-} & \widetilde{G}_u^{^-}  &  G_u^{^+} & \widetilde{G}_u^{^+} & \gamma_u &\\ 
\cline{1-7}
0 & 2^{10}3^45.7 & 1 & 1 \ & \  E_7 & E_7 & 1 &\\
\cline{1-7}
1 & 2^{10}3^25 & A_1 & A_1&\ D_6& D_6  & 1 & \\
\cline{1-7}
2 & 2^{10}3& (A_1)^2 & B_2 & A_1\times D_4 & A_1\times B_4 & 2 &\\
\cline{1-7}
3 &  2^{10}3^2 & (A_1)^3  &\ \ B_3 \ \  & D_4 \  &  F_4 \   & 3 &\\
\cline{1-7}
\ 3' & 2^83 & (A_1)^3  & B_3  & (A_1)^4 & A_1\times B_3  \  & 3 & \\
\cline{1-7}
4 &  2^{10}3^2 &D_4 \  & F_4 \  & (A_1)^3 & B_3 \  & 3 \\
\cline{1-7}
\ 4' & 2^83 & (A_1)^4 & A_1 \times B_3   \  & (A_1)^3 & B_3 \ & 3 &\\
\cline{1-7}
5& 2^{10}3 &A_1\times D_4 & A_1\times B_4 & (A_1)^2 & B_2 &2 &  \\
\cline{1-7}
6 &2^{10}3^25 & D_6& D_6  & A_1 & A_1 & 1 &\\
\cline{1-7}
7 &2^{10}3^45.7 & \  E_7 & E_7 & 1 & 1 & 1&\\
\cline{1-7}
\end{array}$$
\noi {\it Vérification du tableau.}

\smallskip
 Soit $d = \deg(u)$. Comme $G$ contient $-1$, il nous suffit de traiter les cas où $d \leqslant 3$; les autres s'en déduisent en remplaçant $u$ par $-u$; noter que, si $u$ est de type $3$ (resp. $3'$), $-u$ est de type $4$ (resp. $4'$).
 
 Les cas $d=0$ et $d=1$  sont immédiats. 
 
 Lorsque $d=2$ , les prop.8.1 et 8.2
 entraînent que  $G_u^{^-}$ est de type $(A_1)^2$, $\widetilde{G}_u^{^-}$ est de type $B_2$ et $\Gamma_u \simeq \Sym_2$. Le groupe $G_u^{^+}$ est de type $A_1\times D_4$ d'après le §3; comme  $ \widetilde{G}_u^{^-}$ contient $G_u^{^+}$ avec indice 2, il est de type $A_1\times B_4$.

Supposons que $d=3$. Comme ci-dessus, les prop.8.1 et 8.2 entraînent que  $G_u^{^-}$ est de type $(A_1)^3$, $\widetilde{G}_u^{^-}$ est de type $B_3$ et $\Gamma_u \simeq \Sym_2$. 

 Quant $u$ est de type  3, le même argument que pour $d=2$ montre que  $G_u^{^+}$  est d'ordre $2^63$, donc de type $D_4$. Comme c'est un sous-groupe normal de
 $ \widetilde{G}_u^{^+}$, et que $ \widetilde{G}_u^{^+}/G_u^{^+}$ est isomorphe
 à $\Sym_3$, on en déduit que $ \widetilde{G}_u^{^+}$ est de type $F_4$.
 
 Quand $u$ est de type $3'$, $G_u^{^+}$ est d'ordre $2^4$, donc de type $(A_1)^4$, alors que $ \widetilde{G}_u^{^+}$ est d'ordre $2^53$, et le contient comme
 sous-groupe normal d'indice $6$; cela entraîne que $ \widetilde{G}_u^{^+}$ est de type $A_1 \times B_3$.
 Noter que le facteur $A_1$ est contenu dans le noyau de la surjection $ \widetilde{G}_u^{^+}\to \Sym_3$, donc est engendré par une $G$-réflexion; c'est l'un des facteurs de $G_u^{^+}$.
 
   \bigskip

\noi {\bf 11.  Type $E_8$}.

\smallskip

Pour $d= 0,1,2,3$ et $d=5,6,7,8$ il y a une seule classe d'involutions de degré $d$.
Pour $d=4$, il y en a deux: celle appelée {\it du type rectangle} dans [Se 22], 5.8,
et celle appelée {\it du type tétraèdre}; dans le tableau ci-dessous, elles correspondent
aux lignes 4 et  $4'$.

$$\begin{array}{|c|c|c|c|c|c|c|l}
\cline{1-7}
\hbox{$\deg(u)$} & |G_u| &  G_u^{^-} & \widetilde{G}_u^{^-}  &  G_u^{^+} & \widetilde{G}_u^{^+} &\gamma_u & \\ 
\cline{1-7}
0 & 2^{14}3^55^27 & 1 & 1 \ & \  E_8 & E_8 & 1& \\
\cline{1-7}
1 & 2^{11}3^45.7 & A_1 & A_1& E_7& E_7  & 1& \\
\cline{1-7}
2 & 2^{12}3^25 & (A_1)^2 & B_2 & D_6 & B_6 & 2&\\
\cline{1-7}
3 &  2^{11}3^2 & (A_1)^3  & B_3\  &A_1\times D_4 \  &  A_1 \times F_4 \  & 3 &\\
\cline{1-7}
4 &2^{13}3^3 & D_4 \  & F_4 \  & D_4 \  & F_4 \ & 3 & \\
\cline{1-7}
\ 4' & 2^{11}3 & (A_1)^4 & B_4 \  & (A_1)^4 & B_4 \   & 4 &\\
\cline{1-7}
5 &2^{11}3^2 &A_1\times D_4 \  &  A_1 \times F_4 \   & (A_1)^3  & B_3 & 3 & \\
\cline{1-7}
6 & 2^{12}3^25  &D_6 & B_6 &(A_1)^2 & B_2 & 2 &\\
\cline{1-7}
7 & 2^{11}3^45.7 & E_7& E_7 &A_1 & A_1 &1&\\
\cline{1-7}
8 & 2^{14}3^55^27 &  E_8 & E_8 & 1 & 1 & 1 &\\
\cline{1-7}
\end{array}$$
\noi {\it Vérification du tableau.}

\smallskip
\noi On peut se borner au cas $d \leqslant 4$ puisque $G$ contient~$-1$. La méthode est la même que pour le type $E_7$. Les cas  $d=0$ et $d=1$ sont immédiats.  Lorsque $d=2$ ou $d=3$, les prop.8.1 et 8.2 entraînent que $G_u^{^-}$  est de type $(A_1)^d$, que $\Gamma_u \simeq \Sym_d$, que $\widetilde{G}_u^{^-}$ est de type $B_d$, et que $\Gamma_u$ est engendré par les images des involutions de $G_u$ de degré 2. 

Pour $d=2$, on a vu au §3 que $G_u^{^+}$ est de type $D_6$; comme $\widetilde{G}_u^{^+}$ le contient avec indice 2, il est de type $B_6$. On a $|G_u| = |G_u^{^-}|\cdot |\widetilde{G}_u^{^+}| = 2^2\cdot2^66!= 2^{12}3^35$.

Pour $d=3$, un argument analogue montre que $G_u^{^+}$ est de type $A_1\times D_4$, que $\widetilde{G}_u^{^+}$ est de type $A_1\times F_4$, et que $|G_u|=2^{11}3^2.$

\noi {\small [Noter
que l'inclusion de $E_7$ dans $E_8$ transforme une involution $v$ de type 3 de $H=E_7$ en une involution $u$ de $G$ de degré 3; le groupe $\widetilde{H}_v^{^+}$ est donc contenu dans $\widetilde{G}_u^{^+}$ , ce qui entraîne que $\widetilde{G}_u^{^+}$ contient $F_4$.]}

\smallskip
Pour $d=4$, on a vu au §3 que $G_u^{^+}$ est, soit de type $D_4$, soit de type $(A_1)^4$, et que les deux cas sont possibles. C'est le premier cas que nous avons choisi de noter 4 dans le tableau ci-dessus, le second étant noté $4'$. Dans les deux cas, $u$ et $-u$ sont conjuguées (cela résulte de la description
de ces classes donnée dans [Se 22], 5.6); d'où $G_u^{^-} \simeq G_u^{^+}$. 
      
      Si $u$ est de type $4'$, la prop.8.1 montre  que $\widetilde{G}_u^{^-}$ est de type $B_4$. On a $|G_u| = |G_u^{^-}|\cdot |\widetilde{G}_u^{^-}| = 2^4\cdot |B_4|=2^{11}3$, et $\Gamma_u = \Sym_4$, ce qui justifie la ligne ($4')$ du tableau. De plus,
une comparaison avec le tableau de $E_7$ montre que l'involution $u$ provient,
par l'injection $E_7 \to E_8$, d'une involution de $E_7$ de type $4'$.
        
    Supposons que $u$ est de type $4$. Cette involution provient d'une involution 
    de type 4 de $E_7$. On en déduit que     $\widetilde{G}_u^{^-}$ contient un sous-groupe de type $F_4$, donc que $\Gamma_u$ contient un sous-groupe isomorphe à
    $\Sym_3$. Or la prop.2.7 montre que $\Gamma_u$ est isomorphe à un sous-groupe de $\Out_c(G_u^{^-}) \simeq \Out_c(D_4)\simeq \Sym_3$. On en conclut 
    que  $\widetilde{G}_u^{^-} \simeq F_4$. [Autre démonstration : utiliser le fait que $F_4$ est un sous-groupe fini maximal de $\GL_4(\Q)$, cf. [Da 65], (4.3).] On déduit de là que $|G_u| = 2^{13}3^3$ et que $\Gamma_u \simeq \Sym_3$, ce qui justifie la ligne 4 du tableau.
         De plus, cet argument montre que les éléments d'ordre 2 de $\Gamma_u$ sont des images d'involutions de degré 2 de $G_u$.

\smallskip

  Il reste à prouver que {\it les involutions de type} $4 \ ($resp. $4'$) {\it sont des rectangles} (resp. {\it des tétraèdres)} au sens de [Se 22], 5.6. Il suffit de le faire pour le type 4, et sur un exemple explicite. Avec les notations de la fin du §8, choisissons $u=s_2 s_5s_7s_0$, qui est de type rectangle, on l'a vu.
  Posons $v = s_2 s_5s_7$; c'est une involution de $E_7$ de type droite, cf. §8.
  Soit $e$ l'élément $-1$ de $E_7$; le produit $ev$ est une involution de $E_7$ de type 4. Son image dans $E_8$ est aussi de type 4. Mais l'image de $e$ dans $E_8$ est $-s_0$; celle de  $ev$  est donc égale à $-u$. On en conclut que $-u$ est de type 4, et la même chose est vraie pour $u$, puisque $u$ et $-u$ sont conjuguées.

    \bigskip

\noi {\bf 12.  Type $F_4$}.

\smallskip

Ce groupe contient $-1$. Pour $d=1,2,3$, il y a deux classes de conjugaison d'involutions de degré $d$. 

  Les réflexions de chacune des deux classes (notées $L$ et $C$ : ``longues'' et ``courtes'') engendrent un sous-groupe normal de $G$, qui est de type 
 $D_4$;  le quotient de $G$ par ce sous-groupe est isomorphe à $\Sym_3$. Les trois sous-groupes d'ordre 2 de $\Sym_3$ correspondent à trois sous-groupes d'indice 3 de $G$, qui sont des
$\cc$-sous-groupes de type $B_4$. Les classes $C$ et $L$ sont permutées par un automorphisme d'ordre 2 de $G$ (par exemple celui qui est évident sur le diagramme de Coxeter). Chaque classe a 12 éléments.

Le groupe $G_u^{^+}$ associé à une réflexion $u$ est de type $B_3$; lorsque $u$
  est de type $L$, cela résulte de la prop.3.3, car $u$ est conjuguée de la
  réflexion notée  $s_0$ dans cette proposition; le cas où $u$ est de type $C$ résulte du précédent en appliquant un automorphisme de $G$ qui permute $C$ et $L$.

\smallskip
  Les deux classes d'involutions de degré $2$, notées $2$ et $2'$, ont les propriétés suivantes:

Une involution $u$ de type 2 se décompose de deux
  façons différentes en produit de deux réflexions $s$ et $s'$; dans l'une de ces décompositions, $s$ et $s'$ appartiennent à la classe $C$;
  dans l'autre décomposition,  
  $s$ et $s'$ appartiennent à la classe $L$. Le groupe  $G_u^{^-}$ est de type $B_2$. 
  Le nombre de telles involutions est $2.3^2$.

    Une involution $u$ de type $2'$ s
    se décompose de façon unique en $u=ss'$, où $s$ est une réflexion de type $C$ et $s'$ est une réflexion de type $L$. Le groupe  $G_u^{^-}$ est de type $A_1\times A_1$. Le nombre de ces involutions est $2^23^2$.
      
        \smallskip
 \noi {\small [Ces propriétés se démontrent, soit en utilisant les plongements $D_4 \subset B_4 \subset F_4$, soit en utilisant la construction de $G$ par dédoublement à partir de $\Sym_4$, cf. [Se 22], 6.12.] }
 
  \smallskip
  On a le tableau suivant :
  
  $$\begin{array}{|c|c|c|c|c|c|c|l}
\cline{1-7}
\hbox{$\deg(u)$} &   |G_u| &  G_u^{^-} & \widetilde{G}_u^{^-}  &  G_u^{^+} & \widetilde{G}_u^{^+} &\gamma_u & \\ 
\cline{1-7}
0&2^73^2& 1&1&F_4&F_4&1&\\
\cline{1-7}
\ 1 \ {\rm et} \ 1'& 2^53&A_1&A_1&B_3&B_3&1&\\
\cline{1-7}
2&2^6&B_2&B_2&B_2&B_2&1\\
\cline{1-7}
\ 2'&2^4&A_1\times A_1&A_1\times A_1&A_1\times A_1&A_1\times A_1&1&\\
\cline{1-7}
\ 3 \ {\rm et} \ 3'&2^53&B_3&B_3&A_1&A_1&1&\\
\cline{1-7}
4&2^73^2&F_4&F_4&1&1&1&\\
\cline{1-7}
\end{array}$$
\noi {\it Vérification du tableau.}

\medskip 
On peut se borner au cas $d \leqslant 2$ puisque $G$ contient $-1$.
Les deux cas $d=0,1$ sont immédiats.
 
Lorsque $u$ est de type 2, on a vu que
$G_u^{^-} $ est de type $B_2$. Or $u$ et $-u$ sont conjugués (puisque leurs centralisateurs
ont le même ordre). Le groupe $G_u^{^+}$ est donc aussi de type $B_2$. Comme $\widetilde{G}_u^{^-} = G_u/G_u^{^+}$, l'ordre de ce groupe est $2^6/2^3=2^3$; puisqu' il contient
$G_u^{^-}$, qui est d'ordre $2^3$, il lui est égal. On a donc $\Gamma_u=1$.

Le cas du type $2'$ est analogue: les groupes $G_u^{^-}$ et $G_u^{^+}$ sont de type $A_1 \times A_1$; l'ordre de $\widetilde{G}_u^{^-} = G_u/G_u^{^+}$  est $2^4/2^2=2^2$, d'où $\widetilde{G}_u^{^-} = G_u^{^-}$ et $\Gamma_u=1$.
\bigskip

Cette vérification termine la démonstration des énoncés du §1.

  \bigskip
  
  \bigskip

           \begin{center}
      
        {\sc Références}
        
        \end{center}
\noi [Bo 68]  N. Bourbaki, {\it Groupes et alg\`ebres de Lie, Chap. IV-VI}, Hermann, Paris, 1968; English translation, Springer-Verlag, 2002.

\noi [Da 65] E.C. Dade, {\it The maximal finite groups of $4\times 4$ integral matrices}, Ill. J. Math. {\bf9} (1965), 99-122.

\noi [DPR 13] J.M. Douglass, G. Pfeiffer \&  G. Röhrle, {\it On reflection subgroups of finite Coxeter groups},  Comm. Algebra {\bf41} (2013), 2574–2592.

\noi [FV 05] G. Felder  \&  P. Veselov, {\it Coxeter group actions on the complement of hyperplanes and special involutions}, J. Eur. Math. Soc. {\bf7} (2005), 101-116.

\noi [Ho 80] R.B. Howlett, {\it Normalizers of parabolic subgroups of reflection groups}, J. London Math. Soc. {\bf21} (1980), 62-80.

\noi [Se 16] J-P. Serre, {\it Finite Groups $:$ an Introduction}, International Press, Somerville, 2016; seconde édition corrigée, 2022.

\noi [Se 22] --------, {\it Groupes de Coxeter finis $:$ involutions et cubes}, Ens. Math. {\bf68} (2022), 99-133.

\bigskip


\end{document}